\documentclass{conm-p-l}

\usepackage{amsmath}
\usepackage{amscd}
\usepackage{amssymb}

\newtheorem{theorem}{Theorem}[section]

\theoremstyle{definition}

\newtheorem{definition}[theorem]{Definition}

\newtheorem{proposition}[theorem]{Proposition}
\newtheorem{remark}[theorem]{Remark}

\theoremstyle{remark}

\newcommand{\C}{{\mathbb C}}

\newcommand{\End}{ {\rm End}}

\newcommand{\g}{\gamma}
\newcommand{\G}{\Gamma}

\newcommand{\Gn}{\G_n}

\newcommand{\Fock}{{\mathcal F}_X}
\newcommand{\Hx}{\mathbb H_X}

\newcommand{\la}{\lambda}

\newcommand{\orbsym}{H^*_{\text{CR}}(X^n/S_n)}
\newcommand{\porbsym}{H_{X,n}}

\newcommand{\Rgn}{R(\Gn)}

\newcommand{\Xn}{X^{[n]}}
\newcommand{\vac}{|0\rangle}

\newcommand{\Z}{ {\mathbb Z} }

\numberwithin{equation}{section}

\begin{document}

\title
{Universal rings arising in geometry and group theory}
\author{Weiqiang Wang}
\address{Department of Mathematics\\
University of Virginia  \\Charlottesville, VA 22904.}
\email{ww9c@virginia.edu}
\thanks{The author is supported in part by an NSF Grant}


\begin{abstract}
Various algebraic structures in geometry and group theory have
appeared to be governed by certain universal rings. Examples
include: the cohomology rings of Hilbert schemes of points on
projective surfaces and quasi-projective surfaces; the Chen-Ruan
orbifold cohomology rings of the symmetric products; the class
algebras of wreath products, as well as their associated graded
algebras with respect to a suitable filtration. We review these
examples, and further provide a new elementary construction and
explanation in the case of symmetric products. We in addition show
that the Jucys-Murphy elements can be used to clarify the
Macdonald's isomorphism between the FH-ring for the symmetric
groups and the ring of symmetric functions.
\end{abstract}

\maketitle


\section{Introduction}
\subsection{}

The Hilbert scheme $\Xn$ of $n$ points on a complex
(quasi-)projective surface $X$ has provided a natural setup for an
interplay among geometry, representation theory and combinatorics.
The starting point is the construction of a Heisenberg algebra in
terms of incidence varieties which acts irreducibly on the direct
sum over all $n$ of the cohomology groups $H^*(\Xn)$ with
$\C$-coefficients \cite{Na1} (cf. \cite{Gro} for another
construction and \cite{Got1, VW} for motivations). In a nutshell,
the Heisenberg algebra reveals the geometric structures of the
Hilbert schemes in a way similar to the cycle structures
for the symmetric groups. This has recently led to a new approach
toward the cohomology rings of Hilbert schemes.

Being closely related yet complementary to each other, the results on these cohomology
rings can be roughly divided as follows: the connections with
vertex operators and $\mathcal W$ algebras \cite{Lehn, LQW1,
LQW4}; the ring generators (and relations) for the rings
$H^*(\Xn)$ \cite{LQW1, LQW2, LQW5} (also cf. \cite{ES, Bea, Mar}
for a classical approach for some special surfaces); the
constructions of the universal rings which govern the rings
$H^*(\Xn)$ depending on whether or not $X$ is compact \cite{LQW3,
LQW5}; the intertwining relations with the study of Chen-Ruan
orbifold cohomology rings \cite{CR} of the symmetric products
\cite{LS1, LS2, QW, Ru1, Uri, Kau, Got3} as well as the class
algebras of the wreath products \cite{Wa1, FW, EG, LT, Wa3, Wa4}.
We remark that there has been also remarkable
connections of Hilbert schemes with combinatorics \cite{Hai} and
interesting work on the motives of Hilbert schemes \cite{dCM,
Got2}.

Much of the advances in these directions has been made after an
earlier review \cite{Wa2} was written. The main purpose of the
present paper is to review some of these recent progress
complementary to {\em loc. cit.}, with an emphasis on the
appearance of certain universal rings which govern various
algebraic structures. Also it is interesting to compare with
\cite{QW} where the axiomatic nature of the vertex operator
approach developed for Hilbert schemes was formulated explicitly
and similar framework was established for the Chen-Ruan orbifold
cohomology rings of the symmetric products. The
Sections~\ref{sect:symmprod} and \ref{subsect:JM} contain new
constructions and results.

\subsection{}

The Heisenberg algebra provides a distinguished linear
basis for $H^*(\Xn)$ parameterized by multi-partitions (this basis
will be referred to as the Heisenberg monomial basis). It is
natural to ask how the structure constants of the cup product
among these basis elements (up to some suitable scalings) behave.
It turns out that the answer depends essentially on whether $X$ is
projective or quasi-projective\footnote{The use of the terminology
``quasi-projective" excludes ``projective" in the present paper.}.

When $X$ is projective, a stability result formulated in \cite{LQW3} states
that the modified structure constants with respect to a suitably
chosen set of elements in $H^*(\Xn)$ (which include the Heisenberg
monomial basis elements) can be uniquely taken to be independent
of $n$. This has led to a universal ring (termed as the Hilbert
ring associated to $X$) which governs the cohomology rings
$H^*(\Xn)$ for all $n$. These are explained in
Sect.~\ref{sect:proj}. The Hilbert ring actually depends not on
the projective surface $X$, but on the ring $H^*(X)$ and the
canonical class of $X$. While all these statements are described
in the classical language, the machinery used to establish them
replies heavily on the connections with vertex operators.

The stability for the cohomology rings of $H^*\Xn)$ for $X$ {\rm
quasi-projective} requires a quite different formulation. It is
shown \cite{LQW5} that the structure constants with respect to the
Heisenberg monomial basis of $H^*(\Xn)$ are independent of $n$,
for a large class of quasi-projective surfaces (conjecturally
\cite{Wa4} for all quasi-projective surfaces). This has led to a
universal ring (called the FH-ring in \cite{LQW5}) which governs
the cohomology rings $H^*(\Xn)$ for all $n$. This stability result
can in turn be interpreted as giving rise to a surjective ring
homomorphism from $H^*(\Xn)$ to $H^*(X^{[n-1]})$ for each $n$.
This phenomenon has been rather unusual, since there is no natural
algebraic embedding of $X^{[n-1]}$ into $\Xn$ for $m <n$. In fact,
an incidence variety relating $\Xn$ and $X^{[n-1]}$ induces such a
ring homomorphism. This will be explained in
Sect.~\ref{sect:quasiproj}.
\subsection{}

Recall that the Hilbert-Chow morphism $\pi_n: \Xn \rightarrow
X^n/S_n$ from the Hilbert schemes to the symmetric products is a
crepant resolution of singularities, where $S_n$ is the $n$-th
symmetric group. A general principle \cite{DHVW} states that the
geometry of an orbifold is ``equivalent" to the geometry of its
crepant resolutions. Originally, the word ``geometry" here was
expressed in terms of (orbifold) Betti numbers. Over the years,
the ``equivalence" extends to (orbifold) Hodge numbers, (orbifold)
elliptic genera, motivic intergration, and (orbifold) cohomology
rings etc, (see, for example, \cite{Bat, BL, BKR, CR, DL, Kal,
Reid, VW, WaZ, Zas} and the references therein). In particular,
the Chen-Ruan cohomology ring $\orbsym$ of the symmetric products
has attracted much attention, as this is related to the cohomology
rings of Hilbert schemes and vertex operators \cite{LS1, LS2,
LQW5, QW, Ru1, Uri}. This provides nice examples of testing and
verifying Ruan's conjectures \cite{Ru1, Ru2}. We remark that there
has been some further reformulation and development of the
orbifold cohomology rings \cite{FG, AGV}.

It developed that the structures of the rings $\orbsym$ for a
compact complex manifold $X$ are very similar to the cohomology
rings of Hilbert schemes of points on projective surfaces. A
stable ring introduced in \cite{QW}, which are analogous to the
Hilbert ring mentioned above for Hilbert schemes, governs the
Chen-Ruan cohomology rings $\orbsym$. The proof therein relies on
an axiomatic vertex operator approach which are parallel to the
developments in Hilbert schemes. An additional bonus of such an
approach is to provide a more conceptual proof of Ruan's conjecture
on the ring isomorphism between $H^*(\Xn)$ and $\orbsym$ when
$X$ is projective and has a numerically trivial canonical class
(which was also established in \cite{FG, Uri} using \cite{LS2}).

However, the orbifold side is simpler than the crepant resolution
side, and often group-theoretic techniques are applicable. A nice
example is provided by the Jucys-Murphy elements of
symmetric groups \cite{Juc, Mur}
which feature significantly in the study of the rings $\orbsym$ \cite{QW}. In
Sect.~\ref{sect:symmprod}, we shall provide a new and direct
construction (without using vertex operators) which explains the
stability for the rings $\orbsym$. The idea here is to enlarge
further the noncommutative orbifold cohomology rings \cite{FG}
(also cf. \cite{LS2}) of the symmetric products.

\subsection{}

The Hilbert schemes are also resolutions of singularities
of wreath product orbifolds (which are generalized symmetric
products) \cite{Wa1}, where the wreath product $\Gn :=\G^n \rtimes
S_n$ is the semidirect product of the product group $\G^n$ and
$S_n$. The study of the class algebras $R(\G_n)$ of the wreath
products $\G_n$ (which can be regarded as the Chen-Ruan cohomology
ring of the symmetric product of the orbifold $pt/\G$) reveals a
stability similar to the Hilbert schemes for {\em projective}
surfaces \cite{Wa3}. We can further introduce a filtration on
$R(\Gn)$ and study its associated graded algebras $\mathcal
G_\G(n)$ \cite{Wa4} (also compare \cite{EG, FH, LS1, Vas}). It is
shown \cite{Wa4} that the structure constants of these graded
algebras $\mathcal G_\G(n)$ with respect to the conjugacy classes
(which are the analog of the Heisenberg monomial basis) are
independent of $n$. This motivated and corresponds to the
stability explained above of the cohomology rings of Hilbert
schemes when the surface is {\em quasi}-projective. This allows
one to introduce a universal ring $\mathcal G_\G$
(called the FH-ring) which governs the algebras
$\mathcal G_\G(n)$ for all $n$. Both the Chen-Ruan cohomology ring
of the symmetric product and the class algebra of wreath product
are nontrivial generalizations in two directions of the class
algebra of the symmetric group, yet they  can be further unified
under the common roof of the wreath product orbifolds.

Macdonald \cite{Mac} has earlier constructed an isomorphism between the
FH-ring $\mathcal G$ \cite{FH} of the symmetric groups (i.e. when
$\G=1$) and the ring of symmetric functions.
(Actually, the universal ring in each of the  different setups above
has been shown to be a polynomial ring which is isomorphic to
a tensor product of several copies of the ring of symmetric functions.)
In addition, he describes
explicitly the symmetric functions corresponding to the {\em
stable} conjugacy classes in $\mathcal G$. In
Sect.~\ref{sect:wreath} we will show how the Jucys-Murphy elements
can be used to provide further explicit correspondences under
Macdonald's isomorphism.

A better way of reading the present paper is to read first
Sect.~\ref{sect:symmprod} by setting $X$ to be a point and
Sect.~\ref{sect:wreath} by setting $\G$ to be trivial, where the
results and notations would be significantly simplified. In this
way the formulations of the stability in Sections~\ref{sect:proj}
and \ref{sect:quasiproj} may become more transparent.

{\bf Acknowledgments.} Special thanks go to Wei-Ping Li and Zhenbo
Qin for very fruitful collaborations. I also benefit from
stimulating discussions and remarks from many
people, including M.~de~Cataldo, I.~Dolgachev, L.~G\"ottsche,
M.~Haiman, M.~Lehn, E.~Markman, H.~Nakajima, Y.~Ruan, and
J.-Y.~Thibon.

\section{The Hilbert ring of Hilbert schemes for projective
surfaces} \label{sect:proj}
\subsection{Generalities on Hilbert schemes} \label{sec_general}

Let $X$ be a (quasi-)projective complex surface, and $\Xn$ be the
Hilbert scheme of $n$ points on $X$. An element in the Hilbert
scheme $\Xn$ is represented by a length $n$ $0$-dimensional closed
subscheme of $X$. According to Fogarty, $\Xn$ is smooth. We denote
by $H^*(\Xn)$ the cohomology group/ring of $\Xn$ with complex
coefficient. We denote

 $$\Hx = \bigoplus_{n=0}^\infty  H^*(\Xn). $$
The element $1$ in $H^0(X^{[0]}) = \C$ is called the {\it vacuum
vector} and denoted by $|0\rangle$. A non-degenerate
super-symmetric bilinear form $(, )$ on $\Hx$ is induced from the
standard one on $H^*(\Xn)$ defined by
$\displaystyle{(\alpha,\beta) =\int_{\Xn} \alpha\beta}$ for
$\alpha, \beta\in H^*(\Xn)$. For $\mathfrak f \in \End(\Hx)$ we
denote by $\mathfrak f^\dagger \in \End(\Hx)$  its {\it adjoint}
operator.

The Betti numbers of Hilbert schemes for an arbitrary
(quasi-)projective surface was computed in \cite{Got1}. As a
corollary, one obtains the graded dimension of $\Hx$ to be

\begin{eqnarray*}
 \dim_q \Hx := \sum_{n =0}^\infty \dim H^*(\Xn) q^n
 = \prod_{r =1}^\infty \frac
 {(1+q^r)^{h^{odd}(X)}}
 {(1-q^r)^{h^{ev}(X)}}
\end{eqnarray*}
where $h^{odd}(X)$ and $h^{ev}(X)$ denote the dimensions of the
cohomology group of $X$ of odd and respectively even degrees.
Based on this formula, it is suggested \cite{VW} that the space
$\Hx$ can be identified as the Fock space of a Heisenberg algebra
associated to the lattice $H^*(X,\Z)/tor$. This has been
established firmly in terms of incidence varities \cite{Na1}.
Similar results were obtained in \cite{Gro}.

We recall the construction \cite{Na1, Na2} when $X$ is projective.
For $n > 0$ and $\ell \ge 0$, we define $Q^{[n+\ell,n]} \subset
X^{[n+\ell]} \times X \times \Xn$ to be the closed subset
$$\{ (\xi, x, \eta) \in X^{[n+\ell]} \times X \times \Xn \, | \,
\xi \supset \eta \text{ and } {\rm Supp} (I_\eta/I_\xi) = \{ x \}
\}.$$
The linear operator $\mathfrak a_{-n}(\alpha) \in {\rm End}(\Hx)$
(called the creation operator) with $\alpha \in H^*(X)$ is defined
by

$$
\mathfrak a_{-n} (\alpha)(a) = pr_{1*}([Q^{[m+n,m]}] \cdot
 \rho^*\alpha \cdot pr_2^*a).
 $$
for all $a \in H^*(X^{[m]})$, where $pr_1, \rho, pr_2$ are the
projections of $X^{[m+n]}\times X \times X^{[m]}$ to $X^{[m+n]},
X, X^{[m]}$ respectively. Here and below we use $\g \cdot \beta$
or $\g\beta$ to denote the cup product of two cohomology classes
$\g,\beta$, and omit the Poincar\'e duality used to switch between
a homology class and a cohomology class. The annihilation operator
$\mathfrak a_n (\alpha)$ can be defined to be $\mathfrak
a_n(\alpha) = (-1)^n \cdot \mathfrak a_{-n}(\alpha)^\dagger$. Also
set $\mathfrak a_0(\alpha) =0$. Nakajima's theorem \cite{Na2} says
that

$$ [\mathfrak a_n(\alpha), \mathfrak a_m(\beta)] = - n \cdot
\delta_{n, -m}
 \cdot \int_X(\alpha \beta) \cdot {\rm Id}.
$$
Moreover, the space $\Hx$ is an irreducible representation of the
Heisenberg algebra.

\begin{remark} \rm   \label{rem:compact}
The above construction remain to be valid when $X$ is
quasi-projective after some suitable modification, cf. \cite{Na1}.
That is, we have the creation operators $\mathfrak a_{-n}
(\alpha)$ $(n > 0)$ associated to a cohomology class $\alpha$ of
$X$ as above. On the other hand, the annihilation operators
$\mathfrak a_{n} (\g)$ $(n>0)$ are associated to a cohomology
class $\g$ with compact support of $X$.
\end{remark}

\subsection{The Hilbert ring} \label{sec_hilbert}

Given a $\Z_2$-graded finite set $S =S_0 \sqcup S_1$, we denote by
${\mathcal P}(S)$ the set of partition-valued functions $\rho
=(\rho(c))_{c \in S}$ on $S$ such that for every $c \in S_1$, the
partition $\rho(c)$ is required to be {\it strict} in the sense
that $\rho(c) =(1^{m_1(c)} 2^{m_2(c)} \ldots )$ with $m_r(c) = 0$
or $1$ for all $r \ge 1$.

Now let us take a linear basis $S= S_0 \cup S_1$ of $H^*(X)$ such
that $1_X \in S_0$, $S_0 \subset H^{\rm even}(X)$ and $S_1 \subset
H^{\rm odd}(X)$. If we write a partition-valued function $\rho
=(\rho (c))_{c \in S}$, where $\rho(c) =(r^{m_r(c)})_{r \ge 1}
=(1^{m_1(c)} 2^{m_2(c)} \ldots)$, then we define

\begin{eqnarray*}
\ell(\rho)
  &=& \sum_{c \in S} \ell(\rho (c))
  = \sum_{c\in S, r\geq 1} m_r(c) \\
\Vert \rho \Vert
  &=& \sum_{c \in S} |\rho (c)|
  =\sum_{c\in S, r\geq 1} r \cdot m_r(c) \\
{\mathcal P}_n(S)
 &=& \{\rho \in{\mathcal P} (S)\;|\;\; \Vert \rho \Vert =n\} .
\end{eqnarray*}
Given $\rho=(\rho(c))_{c\in S}=(r^{m_r(c)})_{c\in S,r \ge 1} \in
{\mathcal P}(S)$ and $n \ge 0$, we define
\begin{eqnarray}
  {\mathfrak a}_{- \rho(c)}(c)
  &=& \prod_{r \ge 1}
   {\mathfrak a}_{-r}(c)^{m_r(c)}  \nonumber  \\
    %
 {\mathfrak a}_{\rho}(n)
 &=& {\bf 1}_{-(n-\Vert \rho \Vert)}
  \prod_{c\in S}{\mathfrak a}_{-{\rho}(c)}(c) \cdot \vac \in H^*(\Xn)
  \label{eq:basis}
\end{eqnarray}
where we fix the order of the elements $c \in S_1$ appearing in
$\displaystyle{\prod_{c \in S}}$ once and for all, and ${\bf
1}_{-k}$ ($k \ge 1$) denotes $\mathfrak a_{-1}(1_X)^k /k!$. By
definition, ${\mathfrak a}_{\rho}(n)=0$ for $n < \Vert \rho
\Vert$.

\begin{theorem} \cite{LQW3} \label{th:hilbindep}
Let $X$ be a projective surface, and $S$ be a linear basis of
$H^*(X)$ which contains $1_X$. Let $\rho , \sigma \in \mathcal
P(S)$. Then the cup product ${\mathfrak a}_{\rho}(n) \cdot
{\mathfrak a}_{\sigma}(n)$ in $H^*(\Xn)$ can be expressed uniquely
as a linear combination
\begin{eqnarray} \label{eq_structure}
{\mathfrak a}_{\rho}(n) \cdot {\mathfrak a}_{\sigma}(n) =
\sum_{\nu \in \mathcal P(S)} A_{\rho\sigma}^\nu {\mathfrak
a}_{\nu}(n)
\end{eqnarray}
where $\Vert \nu \Vert \le \Vert \rho \Vert +\Vert \sigma \Vert$
and the structure coefficients $A_{\rho\sigma}^\nu$ are
independent of $n$.
\end{theorem}

Note that the cohomology classes ${\mathfrak a}_{\nu}(n)$ with
$\Vert \nu \Vert \le n$ span $H^*(\Xn)$ but are not linearly
independent for $n >1$.

\begin{remark}
Although it is not manifest in the formulation of the statement,
the proof of Theorem~\ref{th:hilbindep} used the connections
between the cup products in $H^*(\Xn)$ and vertex operators
developed in \cite{Lehn, LQW1} in an essential way.
\end{remark}

\begin{definition}  \rm
The {\it Hilbert ring} associated to a projective surface $X$,
denoted by ${\mathfrak H}_X$, is defined to be the ring with a
linear basis formed by the symbols ${\mathfrak a}_\rho$, $\rho \in
{\mathcal P}(S)$ and with the multiplication defined by
${\mathfrak a}_{\rho} \cdot {\mathfrak a}_{\sigma} = \sum_{\nu}
A_{\rho\sigma}^\nu {\mathfrak a}_{\nu}$ where the structure
constants $A_{\rho\sigma}^\nu$ are from (\ref{eq_structure}).
\end{definition}

Note that the Hilbert ring does not depend on the choice of a
linear basis $S$ of $H^*(X)$ containing $1_X$ since the operator
${\mathfrak a}_n(\alpha)$ depends on the cohomology class $\alpha
\in H^*(X)$ linearly. The Hilbert ring ${\mathfrak H}_X$ captures
all the information of the cohomology rings of $\Xn$ for all $n$,
as we easily recover the relations (\ref{eq_structure}) from the
ring ${\mathfrak H}_X$. These observations can be summarized into
the following.

\begin{theorem}  \cite{LQW3} \label{th_stab}
For a given projective surface $X$, the cohomology rings
$H^*(\Xn)$, $n\ge 1$ give rise to a Hilbert ring ${\mathfrak H}_X$
which completely encodes the cohomology ring structure of
$H^*(\Xn)$ for each $n$.
\end{theorem}

We shall see in Sect.~\ref{sect:symmprod} that analogous results
hold in the symmetric product setup which afford simple proofs and
explanations. 
\section{The FH-ring of Hilbert schemes for quasi-projective
surfaces} \label{sect:quasiproj}
\subsection{$n$-independence of the structure constants}

Let $X$ be a smooth quasi-projective surface. As before, let us
take a linear basis $S= S_0 \cup S_1$ of $H^*(X)$ such that $1_X
\in S_0$, $S_0 \subset H^{\rm even}(X)$ and $S_1 \subset H^{\rm
odd}(X)$.

Fix $n \ge 1$. For a given $\rho \in \mathcal P(S)$, we set
${\mathfrak b}_\rho(n) = 0 \in H^*(\Xn)$ if $n < \| \rho \|+
\ell(\rho(1_X))$. If $n \ge \| \rho \|+ \ell(\rho(1_X))$, we
define $\tilde \rho \in \mathcal P(S)$ by putting
$m_r(\tilde{\rho}(c)) =m_r({\rho}(c))$ for $c \in S - \{ 1_X \}$,
$m_1(\tilde{\rho}(1_X)) =n -\| \rho \|- \ell(\rho(1_X))$, and
$m_r(\tilde{\rho}(1_X)) =m_{r-1} (\rho(1_X))$ for $r \geq 2$. Note
that $\| \tilde \rho \| = n$. We define ${\mathfrak b}_\rho(n) \in
H^*(\Xn)$ by

\begin{eqnarray*}
   {\mathfrak b}_\rho(n)
 &=&
  \frac{{\bf 1}_{-(n-\|\rho\|- \ell(\rho(1_X)))}}{\prod_{r \geq 2}
   (r^{m_r(\tilde{\rho}(1_X))} m_r(\tilde{\rho}(1_X))!)}
   \left ( \prod_{c \in S, r \ge 1 \atop c \ne 1_X \,\text{or } r > 1}
   {\mathfrak a}_{-r}(c)^{m_r(\tilde{\rho}(c))}
   \right ) \vac.  \qquad   \label{b.2}
\end{eqnarray*}
where we fix the order of the elements $c \in S_1$ appearing in
the product $\displaystyle{\prod_{c \in S}}$ once and for all. We
remark that the only part in the definition of ${\mathfrak
b}_\rho(n)$ involving $n$ is the factor ${\bf 1}_{-(n-\|\rho\|-
\ell(\rho(1_X)))}$.

As a corollary to the construction of Heisenberg algebra \cite{Na2},
$H^*(\Xn)$ has a linear basis consisting of the classes
${\mathfrak b}_{\rho}(n)$ where $\rho \in \mathcal P(S)$ and
$\Vert \rho \Vert + \ell(\rho(1_X)) \le n.$ Fix a positive integer
$n$ and $\rho, \sigma \in {\mathcal P}(S)$ satisfying $\Vert \rho
\Vert + \ell(\rho(1_X)) \le n$ and $\Vert \sigma \Vert +
\ell(\sigma(1_X)) \le n$. Then we can write the cup product
${\mathfrak b}_{\rho}(n) \cdot {\mathfrak b}_{\sigma}(n)$ as
\begin{eqnarray} \label{str_constant}
{\mathfrak b}_{\rho}(n) \cdot {\mathfrak b}_{\sigma}(n) =
\sum_{\nu \in {\mathcal P}(S)} B_{\rho \sigma}^{\nu}(n) \,\,
{\mathfrak b}_{\nu}(n)
\end{eqnarray}
with structure constants $B_{\rho \sigma}^{\nu}(n)$, where $\Vert
\nu \Vert + \ell(\nu(1_X)) \le n$.

We shall need the following technical definition.
\begin{definition} \label{S}
A smooth quasi-projective surface $X$ is said to satisfy the {\it
S-property} if it can be embedded in a smooth projective surface
$\overline{X}$ such that the induced ring homomorphism
$H^*(\overline{X}) \to H^*(X)$ is surjective.
\end{definition}

Let $X$ be a smooth quasi-projective surface with the
$S$-property. This class of quasi-projective surfaces is very
large, including the minimal resolution of a simple singularity,
the cotangent bundle of a smooth projective curve, and the surface
obtained from a smooth projective surface by deleting a point,
etc.

\begin{theorem} \cite{LQW5} \label{th:hilbquasi}
Let $X$ be a smooth quasi-projective surface satisfying the
S-property. Then all the structure constants
$B_{\rho\sigma}^{\nu}(n)$ given in (\ref{str_constant}) are
independent of $n$.
\end{theorem}

\begin{remark} \rm
The statement in Theorem~\ref{th:hilbquasi} was conjectured
\cite{Wa4} to be valid for every quasi-projective surface. The
proof of Theorem~\ref{th:hilbquasi} given in \cite{LQW5} used in
an essential way relations with Hilbert schemes of points on {\em
projective} surfaces \cite{LQW3, LQW4}. Some new approach will be
needed toward the general conjecture.
\end{remark}

\begin{remark}
The definitions of the elements $\mathfrak a_\rho(n)$ in
Sect.~\ref{sect:proj} and $\mathfrak b_\rho(n)$ above (with
different parametrization and scaling) are made in order to
formulate two different types of stability depending on whether or
not the surface $X$ is projective. For example, Theorem~\ref{th:hilbquasi} is
invalid if $X$ is projective.
\end{remark}

\subsection{The FH-ring of Hilbert schemes}
\par
${}^{}$

Recall from Remark~\ref{rem:compact} that the annihilation
operators $\mathfrak a_{n} (\g)$ $(n>0)$ are associated to a
cohomology class $\g$ with compact support of $X$, when $X$ is
quasi-projective. Set $[x]_c \in H_c^4(X)$ to be the Poincar\'e
dual of the homology class in $H_0(X)$ represented by a point in
$X$. Then $\mathfrak A = -\mathfrak a_1 ([x]_c)$ is induced from
an incidence variety in $X^{[n-1]} \times \Xn$. The following
theorem is an equivalent formulation of
Theorem~\ref{th:hilbquasi}.

\begin{theorem} \cite{LQW5}
Let $X$ be a smooth quasi-projective surface which satisfies the
S-property. Then the linear map $\mathfrak A : H^*(X^{[n]})
\rightarrow H^*(X^{[n-1]})$ is a surjective ring homomorphism.
More explicitly, it sends ${\mathfrak b}_{\rho}(n)$ to ${\mathfrak
b}_{\rho}(n-1)$ for all $\rho$.
\end{theorem}

\begin{remark}
In the case when $X$ is the affine plane, the above theorem was
also contained in \cite{LS1}. The work of Farahat-Higman is
very relevant here \cite{FH, Wa4} (also cf.
Sect.~\ref{sect:wreath} below).
\end{remark}

\begin{definition}
Let $X$ be a smooth quasi-projective surface satisfying the
S-property. We define the {\it FH-ring $\mathcal G_X$} associated
to $X$ to be the ring with a linear basis given by the symbols
$\mathfrak b_{\rho}$, ${\rho} \in \mathcal P(S)$, with the product
given by
\begin{eqnarray*}
{\mathfrak b}_{\rho} \cdot {\mathfrak b}_{\sigma} = \sum_{\nu \in
\mathcal P(S_)} B_{\rho \, \sigma}^{\nu} \,\, {\mathfrak b}_{\nu}.
\end{eqnarray*}
where the structure constants $B_{\rho \sigma}^{\nu}$ come from
(\ref{str_constant}).
\end{definition}

The above results have led to the following.
\begin{theorem} \cite{LQW5}
For a quasi-projective surface $X$ which satisfies the S-property,
the cohomology rings $H^*(\Xn)$, $n\ge 1$ give rise to the FH-ring
${\mathcal G}_X$ which completely encodes the cohomology ring
structure of $H^*(\Xn)$ for each $n$.
\end{theorem}
\section{The stable ring of the symmetric products $X^n/S_n$ for $X$ compact}
\label{sect:symmprod}
\subsection{Generalities on Chen-Ruan cohomology rings}

Let $M$ be a complex manifold of complex dimension $d$ with a
finite group $G$ action. Introduce

$$M\diamond G =\{(g, x) \in G \times M\mid gx=x\} =\bigsqcup_{g \in G} M^g,$$
with a $G$-action given by $h.(g,x) =(hgh^{-1}, hx)$. As a vector
space, we define $H^*(M,G)$ to be the cohomology group of
$M\diamond G$ with $\C$-coefficient, or equivalently,

$$H^*(M,G) = \bigoplus_{g\in G} H^*(M^g).$$
The space $H^*(M,G)$ has a natural induced $G$ action, which is
denoted by $\text{ad} \, h : H^*(M^g) \rightarrow
H^*(M^{hgh^{-1}})$. As a vector space, the orbifold cohomology
group $H^*_{\text{CR}}(M/G)$  is the $G$-invariant part of
$H^*(M,G)$, which is isomorphic to
 $$\bigoplus_{[g] \in G_*} H^*(M^g /Z(g))$$
where $G_*$ denotes the set of conjugacy classes of $G$ and
$Z(g)$ denotes the centralizer of $g$ in $G$.

For $g\in G$ and $x \in M^g$, write the eigenvalues of the action
of $g$ on the complex tangent space $TM_x$ to be $\mu_k =e^{2\pi i
r_k}$, where $0 \le r_k <1$. The {\em degree shift number} or {\em
age} is the rational number $F^g_x =\sum_{k=1}^d r_k$, cf.
\cite{Zas}. It depends only on the connected component $Z$ which
contains $x$, so we can denote it by $F^g_Z$. Then associated to a
cohomology class in $H^r(Z)$, we assign the corresponding element
in $H^*(M,G)$ (and thus in $H^*_{\text{CR}}(M/G)$) a degree of
$r+2 F^g_Z$.

A graded ring structure on $H^*_{\text{CR}}(M/G)$ was introduced by Chen
and Ruan \cite{CR}. This was subsequently reformulated in
\cite{FG} by introducing a ring structure on $H^*(M,G)$ first and
then passing to $H^*_{\text{CR}}(M/G)$ by restriction. We shall
use $\circ$ to denote this product. The ring has the following
property: $\alpha
\circ \beta $ lies in $H^*(M^{gh})$ for $\alpha \in H^*(M^g)$ and
$\beta \in H^*(M^h)$.

\subsection{Enlarging the ring $H^*(X^n,S_n)$}
\label{subsect:enlarge}

Let $X$ be a compact complex manifold of dimension $d$. Our main
objects here are the Chen-Ruan cohomology ring $\orbsym$. Let us
fix a linear basis $S =S_0 \sqcup S_1$ of $H^*(X)$ such that $S_0
\subset H^{ev}(X)$, $1_X \in S_0$, and $S_1 \subset H^{odd}(X)$.

If $U$ is a finite set, we denote by $S_U$ the symmetric group of
permutations on $U$. We denote by $\underline{n}$ the set $\{1,2,
\ldots,n\}$, so $S_{\underline{n}}$ is just our usual symmetric
group $S_n$. For $U \subset \underline{n}$, we regard $S_U$ as a
subgroup of $S_n$. We define $X^U = \{f: U \rightarrow X\}$ with a
natural action by $S_U$. In particular, $X^{\underline{n}} =X^n$.

Following \cite{IK}, we introduce the semigroup of ``partial
permutations'' $PS_n$ as follows. A {\em partial permutation} of
the set $\underline{n}$ is a pair $(\sigma, U)$ which consists of
a finite subset $U \subset \underline{n}$ and an element $\sigma
\in S_U$. Denote by $PS_n$ the set of all partial permutations of
$\underline{n}$. The finite set $PS_n$ is endowed with a natural
semigroup structure by letting the product of two elements
$(\sigma, U)$ and $(\tau,V)$ in $PS_n$ to be $(\sigma\tau, U \cup
V)$. We denote by $\Z[PS_n]$ the semigroup algebra over $\Z$. The
$S_n$ acts on $PS_n$ and thus on $\Z[PS_n]$ by $\text{ad}\; g:
(\sigma,U) \mapsto (g \sigma g^{-1}, gU)$.

We introduce
\begin{eqnarray} \label{eq:comp}
PH^*(X^n,S_n) = \bigoplus_{(\sigma, U) \in PS_n} H^{*-2
F_\sigma}((X^U)^\sigma )
\end{eqnarray}
and define a product $\bullet$ on $PH^*(X^n,S_n)$ as follows. Let
$(\sigma,U)$ and $(\tau,V)$ be in $PS_n$. Given $\alpha \in H^{*-2
F_\sigma}( (X^U)^\sigma )$ and $\beta \in H^{*-2 F_\tau}(
(X^V)^\tau )$, we may identify $\alpha$ as an element in $H^{*-2
F_\sigma}( (X^{U\cup V})^\sigma )$ via the natural embedding
$$H^{*-2 F_\sigma}( (X^U)^\sigma ) \subset H^{*-2 F_\sigma}(
(X^{U\cup V})^\sigma ) $$
where $\sigma \in S_U$ is regarded as a permutation in $S_{U\cup
V}$ by fixing $V\backslash U\cap V$ pointwise. Similarly, we may
regard $\beta$ as an element in $H^{*-2 F_\tau}( (X^{U\cup
V})^\tau)$. We define the product $\alpha \bullet \beta$ to be
$\alpha \circ \beta \in H^{*-2 F_{\sigma\tau}}( (X^{U\cup
V})^{\sigma\tau} )$. Then the associativity of $H^*(X^n,S_n)$
\cite{CR, FG} (together with the trivial associativity of the
union operation of sets) implies that $PH^*(X^n,S_n)$ thus defined
is an associative algebra.

Note that $\text{ad}\; g: H^* ( (X^U)^\sigma ) \rightarrow  H^* (
(X^{gU})^{g \sigma g^{-1}} )$ for $g \in S_n$ and $(\sigma,U) \in
PS_n$ defines an action of $S_n$ on $PH^*(X^n,S_n)$. The product
structure on $PH^*(X^n,S_n)$ is clearly compatible with the
$S_n$-action. We denote by $(\porbsym, \bullet)$ the
$S_n$-invariant subalgebra of $PH^*(X^n,S_n)$.

\subsection{The stability of the rings $\orbsym$}

By associating a partition of $|U|$ (the cardinality of $U$) to an
element $(\sigma, U)$ in $PS_n$, we have a natural
parametrization of the $S_n$-orbits of $PS_n$ by partitions $\la$
with $\| \la \| \leq n$. Just as $\orbsym$ has a linear basis
given by $\mathcal P_n(S)$ (cf. e.g. \cite{QW}), we see that
$\porbsym$ affords a linear basis parameterized by $\rho \in
\mathcal P(S)$ with $\Vert \rho \Vert \le n$ as follows. Given
$\rho = (\rho(s))_{s\in S}\in \mathcal P(S)$ with  $\Vert \rho
\Vert \le n$, define a partition $\bar{\rho}$ of $\Vert \rho
\Vert$ to be $\cup_{s\in S} \rho(s)$; that is, $\bar{\rho}$ is
obtained from $\rho \in \mathcal P (S)$ by forgetting the indices
$s\in S$ and then rearranging the parts in descending order. Take an
element $(\sigma, U) \in PS_n$ with $|U| =\Vert \rho \Vert$ and
$\sigma$ of cycle type $\bar{\rho}$, and take an $x \in
H^*((X^U)^\sigma) \cong H^*(X)^{\otimes \ell(\bar{\rho})}$ with
$\ell( \rho(s) )$ factors equal to $s$ for each $s \in S$ such that
the cycles of $\sigma$ associated to $s$ correspond to the partition
$\rho(s)$. Denote by
$c_\rho(n) \in \porbsym$ the sum of the $S_n$-orbit of $x$. It is
easy to see that $c_\rho (n)$ does not depend on the choice of
$(\sigma, U)$ and $x$, but only on $\rho$ and $n$. The $c_\rho
(n)$, where $\rho \in \mathcal P(S)$ and $\Vert \rho \Vert \le n$,
form a linear basis for $\porbsym$.

Let us write

$$c_\rho (n) \bullet c_\sigma (n) =
\sum_{\nu \in \mathcal P(S)} C_{\rho\sigma}^\nu (n) \; c_\nu (n)$$
summed over $\nu$ with $\Vert\nu \Vert \leq n$, where
$C_{\rho\sigma}^\nu(n)$ denotes the structure constants of the
ring $H_{X,n}$.

Keeping (\ref{eq:comp}) in mind, we define a linear map (for $m
\leq n$)
$$\theta_{n,m}: PH^*(X^n,S_n) \rightarrow PH^*(X^m,S_m)$$
to be the identity map on the component $H^{*-2
F_\sigma}((X^U)^\sigma )$ associated to $(\sigma, U) \in PS_n$ if
$U \subset \underline{m}$, and $0$ if $U \not\subset
\underline{m}$. It is clear that $\theta_{n,m}$ is a ring
homomorphism and it is surjective. Since the homomorphisms
$\theta_{n,m}$ for all $m \leq n$ are compatible, we can define an
inverse limit
$$PH^*(X^\infty,S_\infty)
:=\lim\limits_{\stackrel{\longleftarrow}{n}} PH^*(X^n,S_n)$$
which has an induced algebra structure and an induced action by
the group $S_\infty = \cup_n S_n$. Denote by $H_{X,\infty}$ the
algebra of ${S_\infty}$-invariants in $PH^*(X^\infty,S_\infty)$.
One can construct a linear basis $c_\rho$ with $\rho \in \mathcal
P(S)$ of $H_{X,\infty}$, in the same way as constructing the basis
$c_\rho(n)$'s for $\porbsym$. Write

\begin{eqnarray} \label{eq:univ}
c_\rho  \bullet c_\sigma  = \sum_{\nu} C_{\rho\sigma}^\nu c_\nu.
\end{eqnarray}
The homomorphisms $\theta_{n,m}$ induces surjective ring
homomorphisms, which will be also denoted by $\theta_{n,m}$, from
$\porbsym$ to $ H^*_{X,m}$. This in turn is
compatible with a surjective ring homomorphism
$\theta_n: H_{X,\infty} \rightarrow \porbsym$
by letting $\theta_n (c_\rho) =c_\rho(n)$ if $\Vert \rho \Vert
\leq n$ and $\theta_n (c_\rho) =0$ otherwise. It follows that
$C_{\rho\sigma}^\nu(n) =C_{\rho\sigma}^\nu$ for $\rho,\sigma, \nu$
such that $\Vert\rho \Vert \leq n, \Vert \sigma\Vert \leq n$ and
$\Vert \nu\Vert \leq n$. We have established the following.

\begin{theorem}  \label{th:constant}
The structure constants $ C_{\rho\sigma}^\nu (n)$ are independent
of $n$.
\end{theorem}

The `forgetful' map $f_n: PH^*(X^n,S_n)\rightarrow H^*(X^n,S_n)$
which, when restricted to the component $H^{*-2
F_\sigma}((X^U)^\sigma )$ for $(\sigma,U)\in PS_n$, is given by
the inclusion map $H^{*-2 F_\sigma}((X^U)^\sigma ) \subset H^{*-2
F_\sigma}((X^n)^\sigma )$ where $\sigma$ in the latter is regarded
as an element of $S_n$ by fixing pointwise
$\underline{n}\backslash U.$ Clearly, $f_n$ is indeed a surjective
$S_n$-equivariant ring homomorphism, and thus induces a surjective
ring homomorphism from $\porbsym \rightarrow \orbsym$ which will
be again denoted by $f_n$. We define $\mathfrak p_\rho (n) =
f_n(c_\rho(n))$. It follows from (\ref{eq:univ}) that

\begin{eqnarray} \label{eq:stable}
\mathfrak p_\rho (n) \circ \mathfrak p_\sigma (n) = \sum_{\nu}
C_{\rho\sigma}^\nu  \mathfrak p_\nu (n).
\end{eqnarray}

\begin{remark}
It is quite straightforward to define $\mathfrak p_\rho (n)$
directly, just as we defined $c_\rho (n)$. We denote
$\Fock =\oplus_{n=0}^\infty \orbsym. $
It is well known that one can construct a Heisenberg algebra
associated to the lattice $H^*(X)/tor$ acting irreducibly on
$\Fock$ (cf. e.g. \cite{QW}). One can show that the $\mathfrak
p_\rho (n)$ in the present paper coincides with the one defined in
\cite{QW}, Sect.~4.4, in terms of the Heisenberg monomials
(up to a scale multiple independent of $n$). The $\mathfrak p_\rho (n)$'s
span $\orbsym$ but are not linearly independent when $n>1$.
\end{remark}

We summarize the above into a  diagram of surjective
ring homomorphisms:

\begin{eqnarray*}
H_{X,\infty}=\lim\limits_{\stackrel{\longleftarrow}{n}}\porbsym
\quad \ldots \ldots \longrightarrow H_{X,n}
 & \stackrel{\theta_{n,n-1}}{\longrightarrow} & H_{X,n-1}  \longrightarrow \ldots \\
  {\downarrow f_n} &&{ f_{n-1} \downarrow}  \\
  \orbsym  && H^*_{\text{CR}}(X^{n-1}/S_{n-1})
\end{eqnarray*}


In this way, the algebra $H_{X,\infty}$ with the multiplication
(\ref{eq:univ}) governs the structures of the rings $\orbsym$. The
identity (\ref{eq:stable}) with structure constants
$C_{\rho\sigma}^\nu$ shows that the algebra $H_{X,\infty}$ can be
identified with the {\it stable ring} ${\mathfrak R}_X$ associated
to a compact complex manifold $X$ introduced in \cite{QW}. Thus,
we have established the following theorem which first appeared in
\cite{QW} with the assumption that the complex dimension of $X$ is
even. The proof in {\em loc. cit.} was highly nontrivial and used
a vertex operator method developed therein for the rings
$\orbsym$.
\begin{theorem}
For a compact complex manifold $X$, the Chen-Ruan cohomology rings
$\orbsym$, $n\ge 1$ give rise to a stable ring $H_{X,\infty}$ which completely
encodes the cohomology ring structure of $\orbsym$ for each $n$.
\end{theorem}

\begin{remark}
In the case when $X$ is a point, our constructions of $H_{X,n}$
etc simplify greatly and specialize to that of \cite{IK} which
were used to explain the stability in the class algebras of
symmetric groups discovered in \cite{KO} (also see \cite{Wa3} for
the extension to the wreath product setup.)
\end{remark}

\begin{remark}
Assume that $H^*(X,\Z)$ has no torsion. We can refine the
definition of $\orbsym$ to define the orbifold cohomology ring of
the symmetric product with integer coefficient, which has no
torsion. If we further choose $S$ to be an integral linear basis
for $H^*(X,\Z)$, then the structure constants
$C_{\rho\sigma}^\nu$'s above can be shown to be integers. Also compare
with \cite{AGV} where the orbifold cohomology rings of a general
orbifold can be defined over integers.
\end{remark}
\section{The FH-ring of wreath products}
\label{sect:wreath}
\subsection{Preliminaries on the wreath products}

Let $\G$ be a finite group, and $\G_*$ be the set of conjugacy
classes of $\G$. We will denote
the identity of $\G$ by $1$ and the identity conjugacy class in $\G$ by $c^0$.
The symmetric group $S_n$ acts on the
product group $\Gamma^n= \Gamma \times \cdots \times \Gamma$ by
permutations: $\sigma (g_1, \cdots, g_n)
  = (g_{\sigma^{ -1} (1)}, \cdots, g_{\sigma^{ -1} (n)}).
$ The wreath product of $\Gamma$ with $S_n$ is defined to be the
semidirect product
$$
 \Gamma_n = \{(g, \sigma) | g=(g_1, \cdots, g_n)\in {\Gamma}^n,
\sigma\in S_n \}
$$
 with the multiplication
 $(g, \sigma)\cdot (h, \tau)=(g \, {\sigma} (h), \sigma \tau).$
%
The $i$-th factor subgroup of the product group $\G^n$ will be
denoted by $\G^{(i)}$. The wreath product $\Gn$ embeds in
$\G_{n+1}$ as the subgroup $\Gn \times 1$, and the union $\G_{\infty}
= \cup_{n \ge 1} \Gn$ carries a natural group structure. When $\G$
is trivial, $\G_\infty$ reduces to $S_\infty  = \cup_{n \ge 1}
S_n$.

The space $R_\Z(\G)$ of the class functions of a finite group $\G$
(which is often called the {\em class algebra} of $\G$) is closed
under the convolution. In this way, $R_\Z(\G)$ is also identified
with the center of the group algebra $\Z[\G]$. We will denote by
$\Rgn$ (resp. $R_\Z(\Gn)$) the {\em class algebra} of complex
(resp. integer) class functions on $\Gn$ endowed with the
convolution product.

The conjugacy classes of $\Gn$ can be described in the following
way (cf. \cite{Mac, Zel}). Let $x=(g, \sigma )\in {\Gamma}_n$,
where $g=(g_1, \ldots , g_n) \in {\Gamma}^n,$ $ \sigma \in S_n$.
The permutation $\sigma $ is written as a product of disjoint
cycles. For each such cycle $y=(i_1 i_2 \cdots i_k)$, the element
$p_y =g_{i_k} g_{i_{k -1}} \cdots g_{i_1} \in \Gamma$ is
determined up to conjugacy in $\Gamma$ by $g$ and $y$, and will be
called the {\em cycle-product} of $x$ corresponding to the cycle
$y$. For any conjugacy class $c \in \G_*$ and each integer $i \geq
1$, the number of $i$-cycles in $\sigma$ whose cycle-product lies
in $c$ will be denoted by $m_i(\rho(c))$, and $\rho (c)$ denotes
the partition $(1^{m_1 (c)} 2^{m_2 (c)} \ldots )$. Then each
element $x=(g, \sigma)\in {\Gamma}_n$ gives rise to a
partition-valued function $\rho =( \rho (c))_{c \in \G_*} \in
{\mathcal P} ( \G_*)$ such that $\Vert \rho \Vert : =\sum_{r, c} r
m_r( \rho(c)) =n$. The $\rho$ is called the {\em type} of $x$. It
is known that any two elements of ${\Gamma}_n$ are conjugate in
${\Gamma}_n$ if and only if they have the same type.

\subsection{The modified types} \label{subsect:modifiedtype}

Let $x$ be an element of $\Gn$ of type $\rho =(\rho(c))_{c\in\G_*}
\in \mathcal P_n(\G_*)$. If we regard it as an element in
$\G_{n+k}$ by the natural inclusion $\Gn \leq\G_{n+k}$, then $x$
has the type $\rho\cup (1^k) \in \mathcal P_{n+k}(\G_*)$, where
$(\rho\cup (1^k)) (c) =\rho(c)$ for $c \neq c^0$ and $(\rho\cup
(1^k)) (c^0) =(\rho(c^0), 1, \ldots, 1) =\rho(c^0)\cup (1^k).$ It
is convenient to introduce the {\em modified type} of $x$ to be
$\widetilde{\rho} \in \mathcal P_{n- r}(\G_*)$, where $r
=\ell(\rho(c^0))$, as follows: $\widetilde{\rho}(c) =\rho(c)$ for
$c \neq c^0$ and $\widetilde{\rho}(c^0) = (\rho_1 -1, \ldots,
\rho_r -1)$ if we write the partition ${\rho}(c^0) = (\rho_1,
\ldots, \rho_r)$. Two elements in $\G_\infty$ are conjugate if and
only if their modified types coincide.

Given $\mu \in \mathcal P(\G_*)$, we denote by $\mathcal K_\mu$
the conjugacy class in $\G_\infty$ of elements whose modified type
is $\mu$.
For each $n\ge 0$ and each $\mu \in \mathcal P(\G_*)$, let $K_\mu
(n)$ be the characteristic function of the conjugacy class in
$\Gn$ whose modified type is $\mu$, i.e. the sum of all $\sigma
\in \Gn \cap \mathcal K_\mu$. The nonzero  $K_\mu (n)$'s form a
$\Z$-basis for $R_\Z(\Gn)$.
If $x \in \mathcal K_\mu$, the {\em degree} $\|x\|$ of
$x$ is defined to be the degree $\| \mu\|$ of its {\em modified}
type.

Given $g \in \G$ and $1 \le i \neq j$, we denote by $(i
\stackrel{g}{\rightarrow} j)$ the element
 $((g_1, g_2,\ldots), (i, j))$ in $\G_\infty$, where $(i,j) \in S_\infty$
is a transposition, $g_j=g, g_i =g^{-1}$, and $g_k =1$ for $k \neq
i,j$. Note that $(i \stackrel{g}{\rightarrow} j) =(j
\stackrel{g^{-1}}{\rightarrow} i)$ and it is of order $2.$ The
elements in $\G_\infty$ of the form $(i \stackrel{g}{\rightarrow}
j)$, where $g$ runs over $\G$ and $(i,j)$ (where $i<j$) runs over
all transpositions of $S_\infty$, form the single conjugacy class
whose cycle-products are all $c^0$ and the partition corresponding
to $c^0$ is $(2,1,1,\ldots)$. Clearly, any element $x$ in
$\G_\infty$ can be written as a product of elements in $\mathcal
K_{(1)_{c^0}}$ and elements of the form $h^{(i)} \in \G^{(i)}$, $i
\ge 1$. Such a product is called a {\em reduced expression}
for $x$ \cite{Wa4} if $x$ cannot be written as a product of fewer
such elements. For a general element in $\G_\infty$, a reduced
expression can be constructed cycle-by-cycle.

The number of
elements appearing in a reduced expression of $x$ is shown in {\em
loc. cit.} to be $\Vert \la \Vert$ for a given element $x \in
\G_\infty$ of modified type $\la$. It follows that $\Vert xy \Vert
\le \Vert x \Vert + \Vert y \Vert$, and $\Vert \cdot \Vert$
defines an algebra filtration on $R_\Z(\Gn)$.
The associated graded ring will be denoted by $\mathcal G_\G (n)$.
\subsection{The FH-ring of wreath products}

Given $\la, \mu \in  \mathcal P(\G_*)$, we write the convolution
product $K_ \la (n) \, K_\mu(n)$ in $R_\Z(\Gn)$ as a linear
combination of $K_\nu(n)$:

$$K_\la (n) K_\mu(n) =\sum_\nu D_{\la\mu}^\nu (n) \; K_\nu(n)$$
where the structure constants $D_{\la\mu}^\nu (n) \in \Z_+$
are zero unless $\Vert \nu \Vert \leq \Vert \la\Vert +\Vert
\mu \Vert$.

\begin{theorem}  \label{th:main} \cite{Wa4}
Let $\la,\mu,\nu \in \mathcal P(\G_*)$. Then,
\begin{enumerate}
 \item
there is a unique polynomial $f_{\la\mu}^\nu (x)$ such that
$f_{\la\mu}^\nu (n) =D_{\la\mu}^\nu (n)$ for all positive integers
$n$.
 \item
the polynomial $f_{\la\mu}^\nu (x)$ is a constant if $\Vert \nu
\Vert = \Vert \la\Vert +\Vert \mu \Vert.$
\end{enumerate}
\end{theorem}

The next theorem is a
reformulation of Theorem~\ref{th:main}~(2).

\begin{theorem}  \label{th:fhring}
The graded $\Z$-ring $\mathcal G_\G (n)$ has a multiplication
given by

$$K_\la (n) \, K_\mu(n)
=\sum_{\Vert \nu \Vert = \Vert \la\Vert +\Vert \mu \Vert}
 a_{\la\mu}^\nu  K_\nu(n)$$
where the structure constants $a_{\la\mu}^\nu $ are nonnegative
integers independent of $n$.
\end{theorem}

The results above specialize to the results of Farahat-Higman
\cite{FH} (also cf. \cite{Mac}) in the symmetric group case. Thus,
we make the following definition.

\begin{definition}
The Farahat-Higman ring (or FH-ring) of the wreath products
associated to $\G$, denoted by $\mathcal G_\G$, is a $\Z$-ring
with a $\Z$-basis $K_{\la}$, where $\la \in \mathcal P(\G_*)$, and
a multiplication (called the FH-product)

\begin{eqnarray} \label{eq:const}
K_\la \,  K_\mu
 =\sum_{\Vert \nu \Vert = \Vert \la\Vert +\Vert \mu \Vert}
 a_{\la\mu}^\nu  K_\nu .
\end{eqnarray}
\end{definition}

There is a natural surjective ring homomorphism $\mathcal G_\G
\rightarrow \mathcal G_\G (n)$ for each $n$ which is compatible
with the surjective ring homomorphism $res_n: \mathcal G_\G (n)
\rightarrow \mathcal G_\G (n-1)$ by restriction. The restriction
map $res_n$ clearly sends $K_\la(n)$ to $K_\la (n-1)$ for all
$\la$. The ring $\mathcal G_\G$ can be regarded as the inverse
limit of the family of rings $\{ \mathcal G_\G (n) \}_{n\ge 1}$
under the restriction maps. Theorem~\ref{th:fhring} is
equivalent to the statement that the restriction maps are ring
homomorphisms. We summarize the sequence of surjective ring
homomorphisms above into the following diagram:

\begin{eqnarray*}
\mathcal G_\G = \lim\limits_{\stackrel{\longleftarrow}{n}}
\mathcal G_\G (n) \quad \ldots \ldots \longrightarrow   \mathcal G_\G(n)
\stackrel{res_n}{\longrightarrow} \mathcal G_\G (n-1)
\longrightarrow \ldots
\end{eqnarray*}

\subsection{Macdonald's interpretation of the FH-ring}

We restrict ourselves to the symmetric group case (i.e. when $\G
=1$) in the remainder of this section. In this case,  we shall omit the subscript
$\G$ from the notations for the rings $\mathcal G_\G$, $\mathcal
G_\G(n)$.

Denote by $\Lambda = \oplus_k \Lambda^k$ the $\Z$-ring of
symmetric functions (in infinitely many variables), and
$\Lambda^k$ denotes the subspace of symmetric functions of degree
$k$. Denote by $e_k, h_k, p_k \in \Lambda^k $ respectively the
$k$-th elementary, complete, and power-sum symmetric functions.
There is a standard bilinear form $\langle -,-\rangle$ on
$\Lambda$ and the Schur functions form an orthnormal $\Z$-basis of
$\Lambda$. Denote by $\la '$ the transpose of a partition $\la$.
An involution $\omega: \Lambda \rightarrow \Lambda$ is defined by
$\omega (s_\la) = s_{\la'}$ for all $\la$, and it can also be
characterized as the ring isomorphism of $\Lambda$ by switching
$e_n$ and $h_n$ for all $n$.

A (less standard) involution on the ring $\Lambda$ is defined as
follows \cite{Mac}. Let

$$u = t +h_1 t^2 +h_2 t^3 + \ldots .$$
Then $t$ can be expressed as a power series in $u$, say

$$ t =u + h_1^* u^2 + h_2^* u^3 +\ldots,$$
with coefficients $h_k^* \in \Lambda^k$ ($k \ge 1$). The ring
homomorphism $\psi : \Lambda \rightarrow \Lambda$ defined by $\psi
(h_k) =h_k^*$, $k \ge 1$, is an involution.
We denote $f^* =\psi (f)$. For example, $h_\la^*
=h_{\la_1}^*h_{\la_2}^* \ldots$ for each partition $\la =(\la_1,
\la_2, \ldots)$, and the $h_\la^*$'s form a $\Z$-basis of
$\Lambda$. We denote by $(g_\la )$ the dual basis, i.e. $\langle
g_\la, h_\mu^* \rangle = \delta_{\la,\mu}.$

The following theorem is due to Macdonald \cite{Mac}, Ex.~25, pp.
132 (also cf. Goulden-Jackson \cite{GJ} for another proof).

\begin{theorem} {\rm (Macdonald)} \label{th:macd}
The linear map $\varphi : \Lambda \rightarrow \mathcal G$ defined
by $\varphi (g_\la) = K_\la$ for all partitions $\la$ is a ring
isomorphism.
\end{theorem}

\subsection{Jucys-Murphy elements and Macdonald's isomorphism}
\label{subsect:JM}

Recall \cite{Juc, Mur} that the Jucys-Murphy elements $\xi_{j;n}$
of the symmetric group $S_n$ are defined to be the sums of
transpositions:
\begin{eqnarray*}
\xi_{j;n} = \sum_{i<j} (i,j), \quad j =1, \ldots, n.
\end{eqnarray*}
In particular, we have $\xi_{1;n} =0.$
When it is clear from the text, we may simply write $\xi_{j;n}$ as
$\xi_j$.
An important property of the JM elements is that any
symmetrization of $\{\xi_1, \ldots, \xi_n\}$ lies in $R_\Z(S_n)$.
It is further known that the $k$-th elementary symmetric function
$e_k(\Xi_n)$ of $\Xi_n = \{\xi_1, \ldots, \xi_n\}$ is equal to the
sum of all permutations in $S_n$ having exactly $(n-k)$ cycles.
For example, $e_1(\Xi_n)$ is exactly the characteristic function
$K_{(1)}(n)$ of the conjugacy classes of transpositions in $S_n$.
We will use $\star$ to denote the product in the FH-ring $\mathcal
G(n)$. We define
$$( -\Xi)^{\star k}(n) = \sum_{j=1}^n (-1)^k
\underbrace{\xi_{j;n} \star \ldots \star \xi_{j;n}}_k. $$

\begin{proposition}   \label{prop:exact}
In the ring $\mathcal G(n)$, we have the equality:

\begin{eqnarray*} \label{eq_cocycle}
 ( -\Xi)^{\star k}(n)
  & = &
   \left\{
      \everymath{\displaystyle}
      \begin{array}{ll}
         (-1)^k K_{(k)}(n)  , & \text{if } k<n  \\
          0, &  \text{if } k \geq n .
      \end{array}
    \right.
  \label{cocy}
\end{eqnarray*}
\end{proposition}

\begin{proof}

Since we have $\Vert x \Vert \leq n-1$ for every element $x$ in
$R_\Z(S_n)$, the statement is clear for $k \ge n$.

Now assume $k<n$. We first note that $(i\; j) \star (i\; j) =0$
for any transposition $(i\;j)$, since $(i\; j)^2 =1$, $\Vert
(i\;j)\Vert =1$, and $\Vert 1 \Vert =0$. It follows from the
definition of the Jucys-Murphy elements $\xi_{j;n}$ that the only
permutations which survive in $\underbrace{\xi_{j;n} \star \ldots
\star \xi_{j;n}}_k$ is the product (under the ordinary group
multiplication) of $k$ distinct transpositions and any such a
product is necessarily a $(k+1)$-cycle. Since $( -\Xi)^{\star
k}(n)$ lies in $R_\Z(S_n)$, it is a multiple of $K_{(k)}(n)$. It
remains to show that the multiple is indeed $1$.

The total number of such $(k+1)$-cycles in $\underbrace{\xi_{j;n}
\star \ldots \star \xi_{j;n}}_k$ is $k! \cdot {j-1 \choose k}$,
and thus the total number of such products in $( -\Xi)^{\star
k}(n)$ is $\sum_{j=2}^n k! \cdot {j-1 \choose k}$. The number of
$(k+1)$-cycles appearing in $K_{(k)}(n)$ is $\frac{n!}{(n-k-1)!
\cdot (k+1)}$. We conclude that $( -\Xi)^{\star k}(n) =K_{(k)}(n)$
from the identity (which can be proved easily by induction on $n$)

$$\sum_{j=2}^n k! \cdot {j-1 \choose k} = \frac{n!}{(n-k-1)! \cdot (k+1)}. $$
\end{proof}
Now, the elements $( -\Xi)^{\star k}(n)$, $n\ge 1$, give rise to
an element $( -\Xi)^{\star k}$ in $\mathcal G$.
Let us denote by $\phi  : \mathcal G \rightarrow \Lambda$ the
isomorphism inverse to $\varphi$.

\begin{proposition}  \label{prop:simple}
We have $\phi (K_{(k)}) = -p_k$, and
 $\omega \phi (( -\Xi)^{\star k}) = p_k$.
\end{proposition}

\begin{proof}
The first part was established in \cite{Mac}, Ex.~25~(b), pp. 133.
The second part now follows from the first part,
Proposition~\ref{prop:exact}, and the fact that $\omega( p_k)
=(-1)^{k-1} p_k$.
\end{proof}

Denote by $e_k(-\Xi_n)$, $h_k(-\Xi_n)$ and $s_{\la}(-\Xi_n)$
respectively the $k$-th elementary, complete, and Schur symmetric
function in the variables $-\xi_1, \ldots, -\xi_n$ (using the
FH-product). We denote by $e_k^\star (-\Xi)$, $h^\star _k(-\Xi)$
and $s^\star _\la (-\Xi)$ the corresponding elements in $\mathcal
G$.

\begin{theorem}
We have the following correspondence under the isomorphisms $\phi$
and $\omega \phi$ from $\mathcal G$ to $\Lambda$:
\begin{enumerate}
 \item
  $\phi (K_\la) = g_\la$.
 \item
  $\omega \phi ( (-\Xi)^{\star k}) = p_k$.
 \item
 $\omega \phi (e_k^\star (-\Xi)) = e_k.$
 \item
 $\omega \phi (h_k^\star (-\Xi)) = h_k.$
 \item
 $\omega \phi (s_\la^\star (-\Xi)) = s_{\la}.$
\end{enumerate}
\end{theorem}

\begin{proof}
Part (1) is Theorem~\ref{th:macd} and part (2) is
Proposition~\ref{prop:simple}. The rest follows from part~(2) and
the definitions of $e_k^\star (-\Xi)$, $h^\star _k(-\Xi)$ and
$s^\star _\la (-\Xi)$.
\end{proof}

\begin{remark}
The Jucys-Murphy elements have been vital in understanding the
rings $\orbsym$ for compact $X$ \cite{QW}, and they are intimately
related to vertex operators, cf. \cite{LT, Wa3, QW}. The natural role of the
Jucys-Murphy elements in $\mathcal G$ indicates that they should
play also an important role in the study of the rings $\orbsym$
for noncompact $X$ even without vertex operators.
\end{remark}

\begin{remark}
For an arbitrary finite group $\G$, it has been shown \cite{Wa4}
that the FH-ring $\mathcal G_\G$ is isomorphic to a tensor product
of several copies (parametrized by $\G_*$) of the ring $\Lambda$ of symmetric
functions. It will be very interesting to generalize Macdonald's
symmetric function interpretation of $\mathcal G$
to the general FH-ring $\mathcal G_\G$.
\end{remark}

\end{document}